\documentclass[twoside, a4paper, 12pt]{article}
\UseRawInputEncoding 
\usepackage{amsfonts}
\usepackage{amsmath, amssymb}
\font\tenmsb=msbm10
\newfam\msbfam
\textfont\msbfam=\tenmsb
\def\Bbb#1{{\fam\msbfam#1}}

\title{On the Newton polytope of a Jacobian pair.}
\date{}
\begin{document}
\author {Leonid Makar-Limanov}
% MPIM preprint
\newcommand{\C}{\Bbb C}
\newcommand{\Q}{\Bbb Q}
\newcommand{\Z}{\Bbb Z}
\newcommand{\R}{\Bbb R}
\newcommand{\de}{\partial}
\newcommand{\J}{\mbox{J}}
\newcommand{\Def}{\mbox{def}}
\newcommand{\ord}{\mbox{ord}}
\newcommand{\gap}{\mbox{gap}}
\newcommand{\Frac}{\mbox{Frac}}
\newcommand{\Der}{\mbox{Der}}
\newcommand{\Char}{\mbox{char}}
\newcommand{\ad}{\mbox{ad}}
\newcommand{\al}{\alpha}
\newcommand{\be}{\beta}
\newcommand{\ga}{\gamma}
\newcommand{\ep}{\epsilon}
\newcommand{\la}{\lambda}
\newcommand{\ka}{\kappa}
\newcommand{\Aut}{\mbox{AUT}}
\newcommand{\Ad}{\mbox{Ad}}
\newcommand{\supp}{\mbox{supp}}

\maketitle
\vskip -0cm
\hskip 3cm{{\it To the memory of Anatoli\"{i} Georgievich Vitushkin}}

\hskip 3cm{{\it  one of the champions of the Jacobian Conjecture}}

\vskip 0.5cm

\begin{abstract} The Newton polytope related to a ``minimal" counterexample to the Jacobian conjecture is introduced and described.
This description allows to obtain a sharper estimate for the geometric degree of the polynomial mapping given by a Jacobian pair
and to give a new proof of the Abhyankar's two characteristic pair case.
\end{abstract}



\textbf{Introduction.}\\

Let us assume that $f, \ g \in \C[x, y]$ (where $\C$ is the field of complex numbers) satisfy $\J(f, g) = {\de f\over \de x}{\de g\over \de y} - {\de f\over \de y}{\de g\over \de x}
 = 1$ and is a counterexample to the JC (Jacobian conjecture claiming that $\C[f,g] = \C[x,y]$, see [K]).
It is known for many years that then there exists an automorphism $\xi$ of
$\C[x,y]$ such that the Newton polygon ${\mathcal{N}}(\xi(f))$ of
$\xi(f)$ contains a vertex $v = (m, n)$  where $n > m > 0$ and is
included in a trapezoid with the vertex $v$, edges parallel to the
$y$ axis and to the bisectrix of the first quadrant adjacent to $v$, and two edges belonging to the coordinate axes
(see [A1], [A2], [AO], [GGV], [H], [J], [L], [MW],[M], [Na1], [Na2], [NN1], [NN2], [Ok]). This was improved quite recently by
Pierrette Cassou-Nogu\`{e}s who showed that ${\mathcal N}(f)$ does
not have an edge parallel to the bisectrix (see [CN] and  [ML1]).

So below we assume that ${\mathcal{N}}(f)$ is included in such a trapezoid with the \textit{leading} vertex $(m, n)$.
We may also assume that ${\mathcal{N}}(f)$ and ${\mathcal{N}}(g)$ contain the
origin as a vertex and are similar (easy consequence of the relation $\J(f, g)
= 1$), that  the coefficients with the leading vertices of $f$ and $g$ are equal
to $1$ (this can be achieved by an appropriate re-scaling of $x, \ y$ and $f, \
g$), that $\deg_y(g) > \deg_y(f)$, and that $\deg_y(f)$ does not divide
$\deg_y(g)$ (otherwise we can replace the pair $f, \ g$ by a ``smaller" pair $f,
\  g - cf^k$).

These are the restrictions on ${\mathcal{N}}(f)$  known at present and it is not clear how to further tighten them by working with ${\mathcal{N}}(f)$ only.
To proceed with this line of research I'll consider an irreducible algebraic dependence of $x, \ f, \
g$ and obtain information about the Newton polytope of this dependence.\\

\textbf{Algebraic dependence of $x, \, f$, and $g$.}\\

We can look at $f, \ g$ as polynomials in one variable $y$ over $\C(x)$.
It is well-known that two polynomials in one variable over a field $K$ are algebraically dependent over $K$ (see [W]).
Therefore $f$ and $g$ are algebraically dependent over $\C(x)$.

We may choose a dependence $P(F, G) = P(x, F, G) \in \C(x)[F, G]$ (i.e. $P(x, f, g) = 0$) such that $\deg_G(P)$ is minimal possible
and hence $P$ is irreducible as an element of $\C(x)[F, G]$, with coefficients in $\C[x]$ (since we can multiply a dependence by the
least common denominator of the coefficients), and assume that these polynomial coefficients do not have a common divisor. \\

\textbf{Connection between $G$ and $y$.}\\

$G$ is an algebraic function of $x$ and $F$ given by $P(x, F, G) = 0$ and $y$ is
an algebraic function of $x$ and $F$ given by $F - f(x, y) = 0$.

\textbf{Lemma on $y$}. $y \in \C(x, f, g)$ and $y \in \C(f(c,y), g(c,y))$ for
any $c \in \C$.\\
\textbf{Proof.} By the L\"{u}roth Theorem $\C(f(c, y), g(c, y)) = \C(r(y))$
where $r$ is a rational function (see [W]). We can replace $r$ by its linear
fractional transformation and assume that $r = {p_1(y)\over p_2(y)}$ where $p_1,
\ p_2 \in \C[y]$ and $\deg(p_1) > \deg(p_2)$. Without loss of generality $p_1, \
p_2$  are relatively prime polynomials. Now, $f(c, y) = {F_1(r)\over F_2(r)}$
for some polynomials $F_1, \ F_2$ where $d_1 = \deg(F_1) > d_2 = \deg(F_2)$ and
$f(c, y) = {F_{1, 0}p_1^{d_1} + \dots + F_{1, d_1}p_2^{d_1}\over (F_{2,
0}p_1^{d_2} + \dots + F_{2, d_2}p_2^{d_2})p_2^{d_1 - d_2}}$. Hence $p_2 = 1$ and
$r$ is a polynomial.
Since $1 = \J(f, g)|_{x = c} \in r'(y)\C[y]$ we should have $r'(y) \in \C$.
Therefore $y \in \C(f(c, y), g(c, y))$ and $y \in \C(x, f, g)$. 
Since $x, \ f$ and $g$ are algebraically dependent we can present $y$ as a polynomial in $g$ (with coefficients in $\C(x,f)$). $\Box$\\

\textbf{Remark.} It is easy to prove that $y \in \C(x, f, g)$ using the Jacobian condition only
(${\de f\over \de y} = {P_g\over P_x},  \ {\de g\over \de y} = {-P_f\over P_x}$ since $P(x, f, g) = 0$, hence $\de\over \de y$ acts on $\C(x, f, g)$),
but this does not imply that $y \in \C(f(c,y), g(c,y))$ for all $c \in \C$. $\Box$\\

There is a one to one correspondence between the roots $y_i$ of $f(x, y) - F$ and $G_i$ of $P(x, F, G)$ in any extension of $\C(x, F)$ which contain
these roots. Indeed, $G_i = g(x, y_i)$ and $y_i = R(G_i)$ where $y = R(G) \in \C(x, F)[G]$.\\

\textbf{Newton polyhedron of a polynomial.}\\

Let $p \in \C[x_1, \dots, x_n]$ be a polynomial in $n$ variables. Represent each monomial of $p$ by a lattice point in $n$-dimensional
space with the coordinate vector equal to the degree vector of this monomial. The convex hull ${\mathcal N}(p)$ of the points so obtained
is called the Newton polyhedron of $p$. We will be using this notion in two-dimensional and three-dimensional cases as Newton polygons and Newton polytopes accordingly.\\

\textbf{Weight degree function.}\\

Define a weight degree function on $\C[x_1, \dots, x_n]$ as follows. First, take weights $w(x_i) = \al_i,$ where $\al_i \in \R$ and
put $w(x_1^{j_1} \dots x_n^{j_n}) = \sum_i\al_ij_i$. For a $p \in \C[x_1, \dots, x_n]$ define support $\supp(p)$ as the collection of all
monomials appearing in $p$ with non-zero coefficients. Then $\deg_w(p) = \max(w(\mu) | \mu \in \supp(p))$. Polynomial $p$ can be
written as $p = \sum p_i$ where $p_i$ are forms homogeneous relative to $\deg_w$. The leading form $p_w$ of $p$ according to $\deg_w$ is the
form of the maximal weight of this presentation.

For a non-zero weight degree function monomials appearing in the support of the
leading form of $p$ correspond to the points of a face $\Phi$ of ${\mathcal
N}(p)$ and if the codimension of $\Phi$ is $i$ there is a
cone of dimension $i$ of the weight degree functions corresponding to $\Phi$.
The leading forms corresponding to these weights are the same and we will use
$p(\Phi)$ to denote them.

The correspondence between faces and weight degree functions is
one to one for the faces of the codimension $1$ if we require that the numbers
$\al_1, \dots, \al_n$ are coprime integers. We will some times refer to this
weight degree function as the function corresponding to the face. \\

\textbf{Roots $y_i$ of $F = f(x, y)$.} \\

Newton introduced the polygon which we call the Newton polygon in order to find a
solution $y$ of $p(x, y) = 0$ in terms of $x$ (see [N]).
Here is the process of obtaining such a solution. Consider an edge $e$ of
${\mathcal N}(p)$ which is not parallel to the $x$ axis and take the weight
which corresponds to $e$. Then the leading form $p(e)$
allows to determine the first summand of the solution as follows. Consider an
equation $p(e) = 0$. Since $p(e)$ is a homogeneous form and $\al = w(x) \neq 0$
solutions of this equation are $y = c_i x^{\be\over\al}$ where $\be = w(y)$ and
$c_i \in \C$. Choose any solution $c_i x^{\be\over\al}$ and replace $p(x, y)$ by
$p_1(x, y) = p(x,c_ix^{\be\over\al} + y)$. Though $p_1$ is not necessarily a polynomial in $x$ we
can define the Newton polygon of $p_1$ in the same way as it was done for the
polynomials; the only difference is that $\supp(p_1)$ may contain monomials
$x^\mu y^\nu$ where $\mu \in \Q$ rather than in $\Z$. Further on we will be
using this kind of Newton polygons and Newton polytopes. The polygon ${\mathcal
N}(p_1)$ contains the \textit{degree} vertex $v$ of $e$, i.e. the vertex with
$y$ coordinate equal to $\deg_y(p_w)$ and an edge $e'$ which is a modification
of $e$ ($e'$ may collapse to $v$).
Take the \textit{order} vertex $v_1$ of $e'$, i.e. the vertex with
$y$ coordinate equal to the order of $p_w$ as a polynomial in $y$ (if $e' = v$ take $v_1 = v$). Use
the edge $e_1$ for which $v_1$ is the degree vertex to determine the next summand and so on.
After possibly a countable number of steps we obtain a vertex $v_\mu$ and the edge $e_\mu$ for
which $v_\mu$ is not the degree vertex, i.e. either $e_\mu$ is horizontal or the degree vertex
of $e_\mu$ has a larger $y$ coordinate than the $y$ coordinate of $v_\mu$. It is possible only
if ${\mathcal N}(p_\mu)$ does not have any vertices on the $x$ axis. Therefore $p_\mu(x, 0) = 0$
and a solution is obtained.

When characteristic is zero the process of constructing a solution is more
straightforward then it may seem from this description. The denominators of
fractional powers of $x$ (if denominators and numerators of these rational
numbers are assumed to be relatively prime) do not exceed $\deg_y(p)$. Indeed,
for any initial weight there are at most $\deg_y(p)$ solutions while a summand
$cx^{M\over N}$ can be replaced by $c\varepsilon^Mx^{M\over N}$ where
$\varepsilon^{N} = 1$ which gives at least
$N$ different solutions.

If $\deg_y(p) = n$ and we want to obtain all $n$ solutions we should choose the first edge $e$ appropriately.
Consider $p_w$ where $w(x) = 0, \ w(y) = 1$. This leading form correspond to a horizontal edge with the ``left''
and ``right'' vertices $v_l$ and $v_r$ or a vertex $v$ in case $v_l = v_r$. If we choose $e$ with the degree
vertex $v_r$ we will obtain $n$ solutions with decreasing powers of $x$ and if we choose $e$ with the degree
vertex $v_l$ we will obtain $n$ solutions with increasing powers of $x$. When $v_l = v_r = v$ choose the ``right''
edge containing $v$ to obtain $n$ solutions with decreasing powers of $x$ and 
the ``left'' edge containing $v$ to
obtain $n$ solutions with increasing powers of $x$.

We can apply Newton approach to finding solutions for $F - f(x,y) = 0$ in an
appropriate extension of $\C(x, F)$. To do this we have to take the weights
$w(x), \ w(F), \ w(y)$ so that the corresponding face (possibly an edge) of
${\mathcal{N}}(F - f(x, y))$ contains the leading vertex $(m, n)$ of
${\mathcal{N}}(f(x, y))$ and proceed as above. Of course the process would be
much harder to visualize but
it can be made two-dimensional if the weights $\al = w(x), \ \rho = w(F)$ are
commensurable.
Say, if $w(x) = 0$ replace $\C$ by an algebraic closure $K$ of $\C(x)$ and make
computations over $K$.
If $w(x)\neq 0$ take for $K$ an algebraic closure of
$\C(z)$ where $z = x^{-\rho\over\al}F$, introduce $t$ so that $x = t^{d_1}$ and $F = zt^{-d_2}$ where $d_1, \ d_2 \in \Z$ and ${\al\over\rho} = -{d_1\over d_2}$,
and consider $F - f(x,y) = zt^{-d_2} -
f(t^{d_1}, y)$ as a polynomial in $y, \ t, \ t^{-1}$ over $K$. \\

\textbf{Newton polytope ${\mathcal{N}}(P)$.}\\

In this section we will find some restrictions on ${\mathcal{N}}(P)$.

Observe that  $\deg_y(g^{\deg_y(f)} - f^{\deg_y(g)}) < \deg_y(f)\deg_y(g)$ because of the shape of ${\mathcal{N}}(f)$ and
${\mathcal{N}}(g)$. It is known that the leading form of $P(x, F, G)$ relative to the weight
$w(x) = 0, \ w(F) = \deg_y(f), \ w(G) = \deg_y(g)$ is $p_0(x)(G^{a_0} - F^{b_0})^{\nu}$ where ${a_0\over b_0} = {\deg_y(f)\over\deg_y(g)}, \ (a_0, b_0) = 1$ and
$b_0\nu = \deg_F(P), \ a_0\nu = \deg_G(P)$ (see [PR], [R], [ML2]).

It follows from Lemma on $y$ that $\deg_G(P) = [\C(x, f, g):\C(x, f)] = [\C(x, y):\C(x, f)] = \deg_y(f)$  and that $\deg_G(P_{\lambda}) = \deg_y(f(\lambda, y))$
where $P_{\lambda}$ is an irreducible dependence between $f(\lambda, y)$ and $g(\lambda, y)$ for $\lambda \in \C$ (recall that $y \in \C(x, f, g)$ and  $y \in \C(f(\lambda,y), g(\lambda,y))$.

Furthermore, $\deg_G(P) = \deg_G(P_{\lambda})$ for all $\lambda \in \C^*$ since $\deg_y(f(\lambda, y)) = \deg_y(f)$ for all $\lambda \in \C^*$.
Hence $P_{\lambda}(F, G)$ is proportional to $P(\lambda, F, G)$ for all $\lambda \in \C^*$ and $p_0(\lambda) = 0$ is possible only if $\lambda = 0$.
Therefore $p_0(x) = c_0x^d$ and $(c_0x^d)^{-1}P$ is a polynomial monic in $G$ (with coefficients in $\C[x, x^{-1}]$). From now on $P$ is this monic polynomial.

Denote by $\mathcal{E}$ the edge of ${\mathcal{N}}(P)$ which corresponds to the leading form $(G^{a_0} - F^{b_0})^{\nu}$ of $P$.
This edge belongs to two faces $\Phi_a$ and $\Phi_b$ of ${\mathcal{N}}(P)$. Let us decide that the plane $FOG$ (where $O$ is the origin of the $FGx$ system of coordinates) is horizontal, 
the $x$ axis is vertical, and that $\Phi_a$ is above $\Phi_b$.

The face $\Phi_b$ will be below the plane $FOG$ if $P(x, F, G)$ is a Laurent polynomial in $x$.

Since the leading form of $P$ relative to the weight $w(x) = 0, \ w(F) = \deg_y(f), \ w(G) = \deg_y(g)$ is $(G^{a_0} - F^{b_0})^{\nu}$, the $x$ axis cannot be parallel to $\Phi_a$ or $\Phi_b$.

One can use ${\mathcal{N}}(P)$ to find a presentation of $G$ as a fractional power series in $x, \ F$ using approach discussed in \textbf{Roots $y_i$ of $F = f(x, y)$.} \\

\underline{The face $\Phi_b$.}\\

Assume that the face $\Phi_b$ (the lower face containing $\mathcal{E}$) is below the plane $FOG$. Since the $x$ axis is not parallel to the face $\Phi_b$
we can choose the corresponding weight by taking $w(x) = 1, \ w(F) = \rho < 0, \ w(G) = \sigma  < 0$. Of course, $\rho, \ \sigma \in \Q$. Expansions of $G$
as well as the corresponding expansions of $y$ relative to this weight are by components with the increasing weight.

Consider the leading form $P(\Phi_b)$ and its factorization into irreducible factors. If all these factors depend only on two variables then
$P(\Phi_b) = \phi_1(x, F)\phi_2(x, G)\phi_3(F,G)$ and $\Phi_b$ is either an interval, or a parallelogram, or a hexagon with parallel opposite
sides. Since $\Phi_b$ is neither ($\Phi_b$ is not $\mathcal{E}$ and it cannot contain an edge parallel to $\mathcal{E}$ and of the same length), $P(\Phi_b)$ has an
irreducible factor $Q(x, F, G)$
which depends on $x, \ F$, and $G$. Denote by $\overline{G}$ a root of $Q(x, F, G) = 0$ and by $\widetilde{G}$ a root of $P(x, F, G) = 0$ for
which $\overline{G}$ is the leading form.
Then $f(x, \widetilde{y}) = F$ and $g(x, \widetilde{y}) = \widetilde{G}$ for $\widetilde{y} = R(x, F)[\widetilde{G}]$.
(The reader should think about the roots $\overline{G}, \ \widetilde{G}$, and $\widetilde{y}$ as fractional power series in $x$ over the algebraic closure of the field $\C(z), \ z = x^{-\rho}F$.)

We can write $\widetilde{y} = \sum_{j=0}^\infty y_j$ where $y_j$ are the homogeneous components of $\widetilde{y}$.
Since $f(x, \widetilde{y}) = F$ there exists a $k$ for which $y_j = c_jx^{\mu_j}, \ c_j \in \C, \
\mu_j \in \Q$ if $j \leq k$ and $y_{k + 1} \not\in\overline{\C(x)}$.

We also can get $\widetilde{y}$ from the Newton polytope of $F - f(x, y)$. The terms $y_j$ for $j \leq k$ are obtained by a resolution process applied to ${\mathcal{N}}(f)$ and the term $y_{k + 1}$ is
defined by a face $\Psi$ of this polytope which contains $(0,0,1)$, i.e. the
vertex corresponding to $F$ (otherwise $y_{k + 1} \in \overline{\C(x)}$). The face $\Psi$
corresponds to the weight $w(x) = 1, \ w(F) = \rho, \ w(y) = \alpha = w(y_{k+1})$ and $\Psi$ contains an edge $e \in xOy$ of ${\mathcal{N}}(f(x, \sum_{j=0}^k y_j + y))$.

Denote $f_k(x,y) = f(x, \sum_{j=0}^k y_j + y), \ g_k(x,y) = g(x, \sum_{j=0}^k y_j + y)$
(then ${\mathcal{N}}(f_k)$ contains the edge $e$ and $w(f_k) = \rho$) and
by $f_k(e), \ g_k(e)$ the leading forms of $f_k$ and $g_k$ for the weight $w$. Thus $f_k(e)(x, y_{k+1}) = F$ by definition of $y_{k+1}$;
also $g_k(e)(x, y_{k+1}) \neq 0$ (recall that $y_{k + 1} \not\in \overline{\C(x)}$). Since $g_k(x,\sum_{j = k+1}^\infty y_j) = \widetilde{G}$ we should have $g_k(e)(x, y_{k+1}) = \overline{G}$.

If $\J(f_k(e), g_k(e)) = 0$ then $g_k(e)(x, y_{k+1}) = cF^{\la}, \ c \in \C^*$ (since $f_k(e)$ is a homogeneous form of a non-zero weight any homogeneous form which is
algebraically dependent with $f_k(e)$ is proportional to a rational power of $f_k(e)$). But $\overline{G}$ depends on $x$ and so $\J(f_k(e), g_k(e)) \neq 0$. In view of $\J(f_k, g_k) = 1$ this implies $\J(f_k(e), g_k(e)) = 1$.

Since the expansion $\widetilde{y}$ is by components with the increasing weight, $w(x) > 0, \ w(f_k) < 0$, the leading vertex $(m, n)$ should
be below the line containing $e$. The following consideration shows that this is impossible. We have  $w(g_k) = w(G) = \sigma < 0$ and
$\rho + \sigma = w(x) + w(y)$ to make $\J(f_k(e), g_k(e)) = 1$ possible. Therefore $\rho = w(x) + w(y) - \sigma = 1 + \al - \sigma$ and
points $(\rho, 0)$ and $(1 - \sigma, 1)$ have the same weight $\rho$. (Recall that $w(x) = 1, \ w(y) = \al, \ w(F) = \rho, \ w(G) = \sigma$.)
Thus they both belong to the line containing the edge $e$. But this line intersects the bisectrix of the first quadrant in a point with coordinates
smaller than $1$ since $\rho < 0, \ \sigma < 0$, and the vertex $(m, n)$ is above this line.

Hence $\Phi_b$ cannot be below $FOG$ and $P(x, F, G) \in \C[x, F, G]$. On the other hand $P(0, f(x, 0), g(x, 0)) = 0$ and the Newton polygon of
this dependence is not an edge.{\footnote{If $g(x, 0)^b = cf(x, 0)^a$ then $(f,g)$ cannot be a counterexample to JC because $\C(f(x, 0), g(x, 0)) = \C(x)$.}}  
Therefore $\Phi_b$ is not an edge and belongs to $FOG$.\\

\newpage

\underline{The face $\Phi_a$.}\\

For the face $\Phi_a$, another face which contains $\mathcal{E}$, choose the weight $w(x) = 1, \ w(F) = \rho > 0, \ w(G) = \sigma > 0$.
An expansion of $G$ relative to this weight is by components with the decreasing weight.

Repeating verbatim considerations from the previous subsection we obtain an edge $e$ of the corresponding ${\mathcal{N}}(f_l)$ which
belongs to the line containing points $(\rho, 0), \ (1 - \sigma, 1)$ and running below the leading vertex $(m, n)$.

Therefore $\rho + n[1 - \sigma - \rho] \geq m$, i. e. $n - m \geq n(\rho + \sigma) - \rho$.
Also $\sigma = {b_0\over a_0}\rho$ because $\Phi_a$ contains ${\mathcal{E}}$ and $n - m \geq [n(1 + { b_0\over a_0} ) - 1]\rho$.
Hence $\rho \leq {(n - m)a_0\over n(a_0 + b_0) - a_0}, \ \sigma \leq {(n - m)b_0\over n(a_0 + b_0) - a_0}$ and
$\deg_x(P) \leq n\sigma \leq (n - m) {nb_0\over n(a_0 + b_0) - a_0}$.

If these inequalities are not strict then the edge $e$ contains $(m, n)$ i.e. $e$ is the (right) \textit{leading} edge. Since $\rho < 1, \ \sigma < 1$
this would imply that $f(x, 0)$ and $g(x, 0)$ are constants and then $\J(f, g) = 1$ is impossible.
Therefore $(m, n)$ does not belong to $e$ and the inequalities are strict.

From Lemma on $y$ we have $\C(x, f, g) = \C(x, y)$. Therefore the degree $[\C(x, y):\C(f, g)]$ of the field extension
is equal to $\deg_x(P)$ and 
$$[\C(x, y):\C(f, g)] < (n - m) {nb_0\over n(a_0 + b_0) - a_0}.$$
This estimate is sharper than the estimate $m + n$ obtained by Yitang Zhang (see [Zh]).

It is known that $[\C(x, y):\C(f, g)]$ is at least $6$ if $\J(f,g) = 1$ (see [D1], [D2], [DO], [Or], [S],
[Zo]). Hence the difference $n - m > 6$.\\

\underline{Edges of ${\mathcal{N}}(P)$.}\\

An edge of ${\mathcal{N}}(P)$ can be parallel to a coordinate plane $GOx$ or $FOG$ and then the leading form of $G$ which corresponds to this
edge is $cx^r$ or $cF^r$ where $c \in \C^*, \ r \in \Q$. An edge parallel to $FOx$ does not correspond to a leading form of $G$.

If $E$ is a \textit{slanted} edge i.e. an edge which is not parallel to any coordinate plane then at least one of the corresponding leading forms is $\overline{G} = c x^{r_1} F^{r_2}$
where $c \in \C^*, \ r_i \in \Q^*$. In this case we have more freedom in choosing a weight which corresponds to $E$ and with an appropriate choice the
edge $e \in {\mathcal{N}}(f_k)$ (see \underline{The face $\Phi_b$}) collapses to a vertex and both $f_k(e), \ g_k(e)$ are monomials. Since
$\J(f_k(e), g_k(e)) = 1$ and $\deg_y(f_k(e)), \   \deg_y(g_k(e))$ are non-negative integers either $\deg_y(g_k(e)) = 0$ or $\deg_y(f_k(e)) = 0$.
If $\deg_y(g_k(e)) = 0$ then $\overline{G} = g_k(e)(x, y_{k+1}) = cx^{r_1}$ and the edge $E$ is parallel to $GOx$ and not slanted;
if $\deg_y(f_k(e)) = 0$ then $f_k(e) = cx^s$ while $f_k(e)(x, y_{k+1}) = F$.

Hence ${\mathcal{N}}(P)$ does not have slanted edges. \\

\underline{Non-vertical and non-horizontal faces.}\\

Consider again the face $\Phi_a$. This face belongs to a slanted plane containing $\mathcal{E}$ which intersects the first octant by a triangle
$\triangle$. Since all edges of $\Phi_a$ are parallel to the coordinate planes and $\Phi_a$ contains $\mathcal{E}$, the face $\Phi_a$ is either
$\triangle$ or a trapezoid obtained from $\triangle$ by cutting it with an edge ${\mathcal{E}}_1$ parallel to $\mathcal{E}$.

If $\Phi_a$ is a trapezoid then the same consideration applied to ${\mathcal{E}}_1$ shows that the next face is also a triangle or a trapezoid,
and so on until we reach the face parallel to $FOG$. \\

\underline{Horizontal faces.}\\

The polytope ${\mathcal{N}}(P)$ has a non-degenerate horizontal face $\Phi_b \subset FOG$ (``floor"). It also has a ``ceiling" which may degenerate into a vertex.
Let us replace $f, \ g$ by $f - c_1, \ g - c_2$ where $c_i \in \C$ and $(c_1, c_2)$ is a ``general pair''. Then the corresponding Newton polytope
has a triangular floor (with a vertex in the origin) and a triangular ceiling (with a vertex on the $x$ axis).\\

\underline{The shape of ${\mathcal{N}}(P)$.}\\

Collecting information we obtained about $\mathcal{N}(P)$ we can conclude that all its vertices are in the coordinate planes $FOx$ and $GOx$,
there are two horizontal faces which are right triangles with right angles in the origin and on the $x$ axis, a face $\Phi_G$ in $FOx$ and a face $\Phi_F$
in $GOx$, which are polygons with the same number of vertices, and all remaining faces are trapezoids obtained by connecting the corresponding vertices of
$\Phi_F$ and $\Phi_G$ by edges which are parallel to $\mathcal{E}$.\\

To give a new proof that in the case of two characteristic pairs counterexample is impossible (see [A2]) we will estimate $\rho$ from below.\\

\textbf{An estimate of $\rho$ from below.} \\

In order to get an estimate for $\rho$ of the face $\Phi_a$ from below we should know more about $P(x, F, G)$.

Consider $f, \ g \in \C(x)[y]$. The first necessary ingredient is the expansion of $g$ as a power series of $f$ in an appropriate algebra relative to the weight
given by $w(y) = 1, \ w(x) = 0$.\\

\underline{Expansion of $g$.} \\

Consider the ring $L = \C[x^{-1}, x]$ of Laurent polynomials in $x$. Define $A$ to be the algebra of asymptotic power series in $y$ with coefficients in $L$, i.e.
the elements of $A$ are $\sum^{i=k}_{-\infty} y_i y^i$ where $y_i \in L, \ y_k \ne 0$. For $a = \sum^{i=k}_{-\infty} y_i y^i$ define $|a| = y_k y^k$. \\

\textbf{Lemma on radical.} If $r \in \Q$ is a rational number, $|a| = cx^ly^k, \ c \in \C$, and $|a|^r \in A$ then $a^r \in A$.\\
\textbf{Proof.} From the Newton binomial theorem $a^r = |a|^r\sum_{j=0}^{\infty} {r\choose j}(\sum^{i=k-1}_{-\infty} {y_i\over y_k} y^{i-k})^j$ 
because $a = |a|(1 + \sum^{i=k-1}_{-\infty} {y_i\over y_k} y^{i-k})$.  Since all ${y_i\over y_k} \in L$, element $a^r \in A$. $\Box$\\

Consider $f(x, y), \ g(x, y)$ as elements of $A$. Then $|f| = x^my^n$ and $|g| = c_0|f|^{\la_0}$ where $\la_0 = {b_0\over a_0}$ (see \textbf{Introduction} and
\textbf{Newton polytope ${\mathcal{N}}(P)$}). By lemma on radical $f^{\la_0} \in A$ and hence $g_1 = g - c_0 f^{\la_0} \in A$ (here $c_0 = 1$). Since $\J(f, g_1) = 1$ either
$\J(|f|, |g_1|) = 0$ or $\J(|f|, |g_1|) = 1$. If $\J(|f|, |g_1|) = 0$ then $|g_1| = c_1 |f|^{\la_1}, \ c_1 \in \C, \ \la_1 \in \Q$ and we can define
$g_2 = g - c_0f^{\la_0} - c_1f^{\la_1}$ which is in $A$ for the same reasons as $g_1$. We can proceed until we obtain $g_\ka = g - \sum_{i=0}^{\ka-1} c_if^{\lambda_i} \in A$
for which $\J(|f|, |g_\ka|) = 1$, i.e. $\J(x^my^n, |g_\ka|) = 1$.  Therefore $|g_\ka| = (c_\ka(x^my^n)^{1-n\over n} - {1\over n-m} x^{1-m}y^{1-n})$ where $c_\ka \in \C$.
If $c_\ka \neq 0$ then $(x^my^n)^{1-n\over n} \in A$ and ${m\over n} \in \Z$ which is impossible since $0 < m < n$. Thus $|g_{\kappa}| = {1\over (m-n)} x^{1-m}y^{1-n}$ and 

$$g = \sum_{i=0}^{\ka-1} c_if^{\lambda_i}  + g_\ka, \ c_i \in \C  \ \ \ (1)$$
where $\deg_y(|f^{\la_i}|) > 1-n, \ \deg_y(|g_\ka|) = 1-n$, and $|g_\ka| = {1\over (m-n)} x^{1-m}y^{1-n} = {1\over (m-n)}x^{n-m\over n}|f|^{\la_\ka}$ where $\la_\ka = {1-n\over n}$.

In order to obtain a ``complete" expansion
$$g = \sum_{i = 0}^{\infty} c_if^{\lambda_i} \ \ \  (2)$$
of $g$ through $x$ and $f$ we should extend $A$ to a larger algebra $B$ with elements $\sum^{i=k}_{-\infty} y_i y^i$ where $y_i \in L_n = \C[x^{-m\over n}, x^{m\over n}]$ in which
$f^{1\over n}$ is defined. Indeed $|x^{-m\over n}f^{1\over n}| = y$ and we can obtain an expansion of $g$ with $c_i \in L_n$.

Clearly, $\la_i = {n_i\over n}, \ n_i \in \Z$. Since $\deg_g(P) = n$ and $\la _\ka = {1-n\over n}$ all $n$ roots $G_j$ of $P(x, F, G) = 0$ in $B$ can be obtained from
$G = \sum_{i = 0}^{\infty} c_iF^{n_i\over n}$
by substitutions $F^{1\over n} \rightarrow \varepsilon^j F^{1\over n}, \ j = 0, 1, \dots, n-1$  where $\varepsilon$ is a primitive root of 1 of power $n$.\\

\underline{A monomial of $P(x, F, G)$ containing a power of $x$.}\\

Polytope ${\mathcal{N}}(P)$ contains the edge ${\mathcal{E}}$ with vertices $(n_0, 0,  0)$ and $(0, n, 0)$ where $n_0 = \la_0 n$ (in the system of coordinates $FGx$).
Hence if ${\mathcal{N}}(P)$ contains a vertex $(i, j, k)$ then $\la_0 n \rho \geq i\rho + j\sigma + k = (i + \la_0 j)\rho + k$ and
$\rho \geq {k\over \la_0(n - j) - i}$ which gives a meaningful estimate when $k > 0$.\\

The following algorithm will produce an irreducible relation for polynomials $f, \ g \in \C(x)[y]$.

Put $\tilde{g}_0 = g$. Assume that after $s$ steps we obtained $\tilde{g}_0, \dots, \tilde{g}_s \in \C(x,y)$. Denote $\deg_y(\tilde{g}_i)$ by $m_i$ and the greatest common divisor
of $n, m_0, \dots, m_i$ by $d_i$. Put $d_{-1} = n$ and $a_i = \frac{d_{i-1}}{d_i}$
for $0 \leq i \leq s$. (Clearly $a_sm_s$ is divisible by $d_{s-1}$
and $a_s$ is the smallest integer with this property.) 

Call a monomial ${\bf m} = f^i\tilde{g}_0^{j_0} \dots \tilde{g}_s^{j_s} \ \ s$-\textit{standard} if $0 \leq j_k < a_k, \ k = 0, \dots, s$.
Find an $s-1$-standard monomial ${\bf m}_{s,0}$ with $\deg_y({\bf m}_{s,0}) = a_sm_s$ and $k_0 \in K = \C(x)$ for which $m_{s,1} =
\deg_y(\tilde{g}_s^{a_s} - k_0 {\bf m}_{s,0}) < a_sm_s$.
If $m_{s,1}$ is divisible by $d_s$ find an $s$-standard monomial
${\bf m}_{s,1}$ with $\deg_y({\bf m}_{s,1}) = m_{s,1}$ and $k_1 \in
K$ for which $m_{s,2} = \deg(\tilde{g}_s^{a_s} - k_0 {\bf m}_{s,0} - k_1{\bf m}_{s,1}) < m_{s,1}$ and so on.

If after a finite number of reductions $m_{s, i}$ which is not divisible by $d_s$ is obtained, denote the corresponding
expression by $\tilde{g}_{s+1}$ and make the next step.
After a finite number of steps we obtain an irreducible relation.

This algorithm was suggested in [ML2] with a proof that it works. In the zero characteristic case it is also shown there that all $\tilde{g}_i$
are polynomials in $y$ (i.e. there are no negative powers of $f$ in the standard monomials). 

We can rewrite $(1)$ as
$$g = \sum_{i=0}^{\ka-1} c_if^{n_i\over n}  + g_\ka, \ c_i \in \C  \ \ \ (3)$$
where $|g_{\kappa}| = {1\over (m-n)} {xy\over |f|}$. Applying the algorithm to this expansion we will get after several
steps ``the last'' $\tilde{g}_\ka$ with $|\tilde{g}_\ka| = c|{xy\over f}\tilde{g}_0^{a_0 - 1}\tilde{g}_1^{a_1 - 1}\dots \tilde{g}_{\ka-1}^{a_{\ka - 1} - 1}|$.

In the case of two characteristic pairs $\ka = 1$ and $|\tilde{g}_1| = c|{xy\over f}\tilde{g}_0^{a_0 - 1}|$. If we denote $|f| = (x^ay^b)^{a_0}, \ |g| = (x^ay^b)^{b_0}$
then $P = \tilde{g}_1^b - cx^{b-a} f^i\tilde{g}_0^j - \dots$ where $|x^{b-a} f^i\tilde{g}_0^j| = |{xy\over f}\tilde{g}_0^{a_0 - 1}|^b$.
Therefore $\rho \geq {b - a\over \la_0(n - j) - i} = {b - a\over \la_0(ba_0 - j) - i}$.
Since $|x^{b-a} f^i\tilde{g}_0^j| = |{xy\over f}\tilde{g}_0^{a_0 - 1}|^b = |x^{b-a}(x^ay^b)^{1 - a_0b + b_0(a_0 - 1)b}|$
we have $a_0i + b_0j = 1 - a_0b + b_0(a_0 - 1)b$ and
$i + \la_0 j = {bb_0a_0 - ba_0 - bb_0 + 1\over a_0}$ (recall that $\la_0 = {b_0\over a_0}$).
Hence
$\rho \geq {b - a\over \la_0(ba_0 - j) - i} = {(b - a)a_0\over \la_0 b a_0^2 - (bb_0a_0 - ba_0 - bb_0 + 1)} =
{(b - a)a_0\over b a_0 b_0 - (bb_0a_0 - ba_0 - bb_0 + 1)} = {(b - a)a_0\over ba_0 + bb_0 - 1}$.
On the other hand $\rho < {(n - m)a_0\over n(a_0 + b_0) - a_0} = {(b - a)a_0^2\over ba_0(a_0 + b_0) - a_0} = {(b - a)a_0\over b(a_0 + b_0) - 1}$
and we have a contradiction. \\

\textbf{Acknowledgements.}\\

\begin{small} 
While working on this project the author was supported by the Max-Planck-Institut f\"{u}r Mathematik in Bonn, Germany, NSA and NSF grants,
a Fulbright fellowship awarded by
the United States--Israel Educational Foundation, and a FAPESP grant 2011/52030-5 awarded by the State of S\~{a}o Paulo, Brazil.                                                                                                                                \end{small}
\\

\begin{center}
{\bf References Sited}
\end{center}

\noindent [A1] S. S. Abhyankar, Lectures On Expansion Techniques In Algebraic Geometry, Tata Institute
of Fundamental Research, Bombay, 1977.

\noindent [A2] S. S. Abhyankar, Some remarks on the Jacobian question. With notes by Marius van der Put and William Heinzer. Updated by Avinash Sathaye.
Proc. Indian Acad. Sci. Math. Sci. 104 (1994), no. 3, 515--542.

\noindent [AO] H. Appelgate and H. Onishi, The Jacobian conjecture in two variables, J. Pure Appl. Algebra 37 (1985), no. 3, 215--227.

\noindent [CN] P. Cassou-Nogu\'{e}s, Newton trees at infinity of algebraic curves. Affine algebraic geometry, 1--19, CRM Proc. Lecture Notes, 54, Amer. Math. Soc., Providence, RI, 2011. (The Russell Festschrift.)

\noindent [Di] J. Dixmier, Sur les alg\`{e}bres de Weyl, Bull. Soc. Math. France 96 (1968), 209--242.

\noindent [D1] A. Domrina,  Four-sheeted polynomial mappings in $\C^2$. The general case. (Russian) Mat. Zametki 65 (1999),
no. 3, 464--467; translation in Math. Notes 65 (1999), no. 3-4, 386--389.

\noindent [D2] A. Domrina, Four sheeted polynomial mappings of $\C^2$. II. The general case. (Russian) Izv. Ross. Akad. Nauk Ser. Mat. 64 (2000), no. 1, 3--36; translation in Izv. Math. 64 (2000), no. 1, 1--33.

\noindent [DO] A. Domrina, S. Orevkov, Four-sheeted polynomial mappings of $\C^2$. I. The case of an irreducible ramification curve. (Russian) Mat. Zametki 64 (1998), no. 6, 847--862; translation in Math. Notes 64 (1998), no. 5-6, 732--744 (1999).

\noindent [GGV] J.Guccione, J.Guccione, C.Valqui, On the shape of possible counterexamples to the Jacobian conjecture, J. Algebra, 471 (2017), 13-74.

\noindent [H] R. Heitmann, On the Jacobian conjecture, J. Pure Appl. Algebra 64 (1990), 35--72.

\noindent [J] A. Joseph, The Weyl algebra — semisimple and nilpotent elements. Amer. J. Math. 97 (1975), no. 3, 597--615.

\noindent [K] O.H. Keller, Ganze Cremona-Transformationen, Monatsh. Math. Physik 47 (1939) 299--306.

\noindent [L] J. Lang, Jacobian pairs II, J. Pure Appl. Algebra 74 (1991), 61--71.

\noindent [M] T. Moh, On the Jacobian conjecture and the configuration of roots, J. reine angew. Math., 340 (1983), 140--212.

\noindent [ML1] L. Makar-Limanov, On the Newton polygon of a Jacobian mate. Automorphisms in birational and affine geometry, 469–476, Springer Proc. Math. Stat., 79, Springer, Cham, 2014. 

\noindent [ML2] L. Makar-Limanov, A new proof of the Abhyankar-Moh-Suzuki theorem via a new algorithm for the polynomial dependence. J. Algebra Appl. 14 (2015), no. 9, 1540001, 12 pp.

\noindent [MW] J. McKay, S. Wang, A note on the Jacobian condition and two points at infinity. Proc. Amer. Math. Soc. 111 (1991), no. 1, 35--43.

\noindent [N] I. Newton, De methodis serierum et fluxionum, in D. T Whiteside (ed.), The Mathematical Papers of Isaac Newton, Cambridge University Press, Cambridge, vol. 3, 1967-1981, 32--353; pages 43--71.

\noindent [Na1] M. Nagata, Two-dimensional Jacobian conjecture. Algebra and topology 1988 (Taejon, 1988), 77-98, Korea Inst. Tech., Taejon, 1988.

\noindent [Na2] M. Nagata, Some remarks on the two-dimensional Jacobian conjecture. Chinese J. Math. 17 (1989), no. 1, 1--7.

\noindent [NN1] A. Nowicki, Y. Nakai, On Appelgate-Onishi's lemmas. J. Pure Appl. Algebra 51 (1988), no. 3, 305--310.

\noindent [NN2] A. Nowicki, Y. Nakai, Correction to: "On Appelgate-Onishi's lemmas'' J. Pure Appl. Algebra 58 (1989), no. 1, 101.

\noindent [Ok] M. Oka, On the boundary obstructions to the Jacobian problem. Kodai Math. J. 6 (1983), no. 3, 419--433.

\noindent [Or] S. Orevkov, On three-sheeted polynomial mappings of $\C^2$. (Russian) Izv. Akad. Nauk SSSR Ser. Ma\t. 50 (1986), no. 6, 1231–1240, 1343.

\noindent [PR] B. Peskin; D. Richman, A method to compute minimal polynomials. SIAM J. Algebraic Discrete Methods 6 (1985), no. 2, 292-299.

\noindent [Ri] D. Richman,  On the computation of minimal polynomials. J. Algebra 103 (1986), no. 1, 1-17.

\noindent [S] I. Sigray, Jacobian trees and their applications. Thesis (Ph.D.), E\"{o}tv\"{o}s   Lor\'{a}nd
University (ELTE), Budapest, Hungary, 2008.

\noindent [W] B. van der Waerden, Algebra. Vol. I. Based in part on lectures by E. Artin and E. Noether. Translated from the seventh German edition
by Fred Blum and John R. Schulenberger. Springer-Verlag, New York, 1991. xiv+265 pp.

\noindent [Zh]  Y. Zhang, The Jacobian conjecture and the degree of field extension. Thesis (Ph.D.), Purdue University. 1991.

\noindent [Zo] H. \.{Z}o\l\c{a}dek, 
An application of Newton-Puiseux charts to the Jacobian problem. 
Topology 47 (2008), no. 6, 431--469.

%\newpage

\vspace{1cm}

\begin{footnotesize}\noindent Department of Mathematics, Wayne State University,
Detroit, MI 48202, USA; \\ Department of Mathematics \& Computer Science, the Weizmann Institute of Science,
Rehovot 76100, Israel.                           \end{footnotesize}
\vskip 0.25cm
\begin{footnotesize}\noindent \textit{E-mail address}: lml@wayne.edu\\                                                       \end{footnotesize}

\end{document}